\documentclass[10pt,a4paper]{article}
\usepackage{amssymb}
\usepackage{latexsym}
\newtheorem{theorem}{Theorem}[section]

\newtheorem{lemma}{Lemma}[section]
\newtheorem{proposition}{Proposition}[section]
\newtheorem{remark}{Remark}[section]

\newtheorem{corollary}{Corollary}[section]
\def\proof{\mbox {\it \textbf{Proof.}~~}}
\usepackage{amsmath}
\numberwithin{equation}{section}
\begin{document}
\title{{\bf\Large  Multiple solutions for a Kirchhoff-type equation with general nonlinearity}}
\author{{\small Sheng-Sen Lu}
\vspace{1mm}\\
{\it\small Center for Applied Mathematics}, {\it\small Tianjin University}\\
{\it\small Tianjin, 300072,}\\
{\it\small and}\\
{\it\small Chern Institute of Mathematics and LPMC}, {\it\small Nankai University}\\
{\it\small Tianjin, 300071, PR China}\\
{\it\small e-mail: lushengsen@mail.nankai.edu.cn}}\vspace{1mm}
\date{}
\maketitle
\begin{center}
{\bf\small Abstract}
\vspace{3mm}
\hspace{.05in}\parbox{4.5in}
{{\small This paper is devoted to the study of the following autonomous Kirchhoff-type equation
$$-M\left(\int_{\mathbb{R}^N}|\nabla{u}|^2\right)\Delta{u}= f(u),~~~~u\in H^1(\mathbb{R}^N),$$
where $M$ is a continuous non-degenerate function and $N\geq2$. Under suitable additional conditions on $M$ and general Berestycki-Lions type assumptions on the nonlinearity $f$, we establish several existence results of multiple solutions by variational methods, which are also naturally interpreted from a non-variational point of view.}}
\end{center}

\section{Introduction and main results}
In this paper, we consider the following autonomous nonlinear elliptic problem with a general subcritical nonlinearity:
\begin{equation}\tag{$\mathcal{KT}$}
\left\{
\begin{aligned}
&-M\left(\int_{\mathbb{R}^N}|\nabla{u}|^2\right)\Delta{u}= f(u)~~\text{in}~\mathbb{R}^N,\\
&u\in H^1(\mathbb{R}^N),~~~~u\not\equiv0~~\text{in}~\mathbb{R}^N,\\
\end{aligned}
\right.
\end{equation}
where $N\geq2$, $M:\mathbb{R_+}\to \mathbb{R}_+$ and $f:\mathbb{R}\to\mathbb{R}$ are continuous functions that satisfy some assumptions which will be specified later on.

In the case where $M$ is not identically equal to a positive constant, the class of Problem $(\mathcal{KT})$ is called of Kirchhoff type because it comes from an important application in Physic and Engineering. Indeed, if we let $M(t)=a+bt$ with $a,b>0$ and replace $\mathbb{R}^N$ and $f(u)$ by a bounded domain $\Omega\subset\mathbb{R}^N$ and $f(x,u)$ respectively in $(\mathcal{KT})$, then we get the following Kirchhoff problem:
\begin{equation*}
-\left(a+b\int_{\Omega}|\nabla{u}|^2\right)\Delta{u}= f(x,u)~~\text{in}~\Omega,
\end{equation*}
assuming the homogeneous Dirichlet boundary condition, which is related to the stationary analogue of the equation
\begin{equation*}
\rho\frac{\partial^2u}{\partial t^2}-\left(\frac{P_0}{h}+\frac{E}{2L}\int^L_0\left|\frac{\partial u}{\partial x}\right|^2\right)\frac{\partial^2u}{\partial x^2}=0~~\text{in}~(0,T)\times(0,L)
\end{equation*}
presented by G. Kirchhoff in \cite{Ki83}. Besides, $(\mathcal{KT})$ is also called a nonlocal problem in this case because of the appearance of the term $M\left(\int_{\mathbb{R}^N}|\nabla{u}|^2\right)\Delta{u}$ which implies that $(\mathcal{KT})$ is no longer a pointwise identity. And this phenomenon provokes some mathematical difficulties which make the study of Problem $(\mathcal{KT})$ particularly interesting.

On the other hand, when $M$ is identically equal to a positive constant, for example $M(t)\equiv 1$, there has been a considerable amount of research on this kind of problems during the past years. The interest comes, essentially, from two reasons: one is the fact that such problems arise naturally in various branches of Mathematical Physics, indeed the solutions of $(\mathcal{KT})$ in the case where $M(t)\equiv 1$ can be seen as solitary waves (stationary states) in nonlinear equations of the Klein-Gordon or Schr\"{o}dinger type, and the other is the lack of compactness, a challenging obstacle to the use of the variational methods in a standard way.

In the celebrated papers \cite{Be83-3,Be83-1,Be83-2}, the authors studied the case where $M(t)\equiv 1$, namely the following autonomous nonlinear scalar field problem
\begin{equation}\tag{$\mathcal{SF}$}
\left\{
\begin{aligned}
&-\Delta{u}= f(u)&\text{in}~\mathbb{R}^N,\\
&u\in H^1(\mathbb{R}^N),~~u\not\equiv0&\text{in}~\mathbb{R}^N,
\end{aligned}
\right.
\end{equation}
under the following assumptions on the nonlinearity $f$:
\begin{itemize}
  \item[$(f_0)$] $f\in C(\mathbb{R},\mathbb{R})$ is continuous and odd.
  \item[$(f_1)$] For $N\geq3$, we have
      \begin{equation}\label{equ1.1}
      -\infty<\underset{t\rightarrow0}{\liminf}\frac{f(t)}{t}\leq\underset{t\rightarrow0}{\limsup}\frac{f(t)}{t}<0.
      \end{equation}
      For $N=2$, we have
      \begin{equation}\label{equ1.2}
      \underset{t\rightarrow0}{\lim}\frac{f(t)}{t}\in(-\infty,0).
      \end{equation}
  \item[$(f_2)$] When $N\geq3$, we have $$\underset{t\rightarrow\infty}{\lim}\frac{f(t)}{|t|^{\frac{N+2}{N-2}}}=0.$$

  When $N=2$, for any $\alpha>0$, we have
      \begin{equation*}
      \underset{t\rightarrow\infty}{\lim}\frac{f(t)}{e^{\alpha t^2}}=0.
      \end{equation*}
  \item[$(f_3)$] There exists $\zeta>0$ such that $F(\zeta)>0$, where $F(t):=\int^t_0f(\tau)d\tau$.
\end{itemize}

With the aid of variational methods and critical points theory, by studying certain constrained problems, Berestycki-Lions and Berestycki-Gallouet-Kavian established the existence results of a ground state, namely a nontrivial solution which minimizes the action functional among all the nontrivial solutions, and infinitely many bound state solutions of $(\mathcal{SF})$ in \cite{Be83-1,Be83-2} for $N\geq3$ and in \cite{Be83-3} for $N=2$ respectively.

As we can see, there is a difference in the assumption $(f_1)$ between the cases $N\geq3$ and $N=2$. We remark here that, in the proofs given by \cite{Be83-3} for the case $N=2$, the existence of a limit $\lim_{t\rightarrow0}f(t)/t\in(-\infty,0)$ is used in an essential way to show the Palais-Smale compactness condition for the corresponding functional under suitable constraints. It is hard to generalize \eqref{equ1.2} to the general inequality \eqref{equ1.1} in that argument.

Later on, in a recent paper \cite{Hi10}, Hirata, Ikoma and Tanaka revisited Problem $(\mathcal{SF})$ in the case $N\geq2$ assuming $(f_0)$, $(f_2)$, $(f_3)$ and
\begin{itemize}
  \item[$(f'_1)$] $-\infty<\underset{t\rightarrow0}{\liminf}\frac{f(t)}{t}\leq\underset{t\rightarrow0}{\limsup}\frac{f(t)}{t}<0$
\end{itemize}
and managed to find radial solutions through the unconstrained functional
\begin{equation}\label{equ1.3}
I(u):=\frac{1}{2}\int_{\mathbb{R}^N}|\nabla u|^2-\int_{\mathbb{R}^N}F(u),~~~~u\in H^1(\mathbb{R}^N).
\end{equation}
In \cite{Hi10}, following the approach introduced by Jeanjean in \cite{Je97}, Hirata, Ikoma and Tanaka considered the auxiliary functional $\tilde{I}:\mathbb{R}\times H^1_r(\mathbb{R}^N)\to\mathbb{R}$
\begin{equation*}
\tilde{I}(\theta,u):=\frac{1}{2}e^{(N-2)\theta}\int_{\mathbb{R}^N}|\nabla u|^2-e^{N\theta}\int_{\mathbb{R}^N}F(u).
\end{equation*}
In this way, they were able to find a Palais-Smale sequence $(\theta_j,u_j)^{+\infty}_{j=1}$ in the augmented space $\mathbb{R}\times H^1_r(\mathbb{R}^N)$ such that $\theta_j\to 0$ and $u_j$ "almost" satisfies the Poho\u{z}aev identity associated to $(\mathcal{SF})$. With the aid of this extra information, it was proved that Problem $(\mathcal{SF})$ possesses a positive least energy solution and infinitely many (possibly sign changing) radially symmetric solutions.

Our main goal of the present paper is to try to provide some multiplicity results for Problem $(\mathcal{KT})$ under the very general assumptions $(f_0)$, $(f_2)$, $(f_3)$ and $(f'_1)$ on $f$ and some suitable conditions on $M$ by variational methods.

In terms of $(f_0)$, $(f'_1)$ and $(f_2)$, we conclude that the corresponding functional $J$ of $(\mathcal{KT})$ given by
\begin{equation*}
J(u):=\frac{1}{2}\widehat{M}\left(\int_{\mathbb{R}^N}|\nabla u|^2\right)-\int_{\mathbb{R}^N}F(u)
\end{equation*}
is well-defined on $H^1(\mathbb{R}^N)$ and of class $C^1$, where $\widehat{M}(t):=\int^t_0M(\tau)d\tau$. It is easy to see that $J$ is invariant under rotations of $\mathbb{R}^N$. Then,
\begin{equation*}
H^1_r(\mathbb{R}^N):=\left\{u\in H^1(\mathbb{R}^N)~|~u(x)=u(|x|)\right\}
\end{equation*}
is a natural constraint to look for critical points of $J$, namely critical points of the functional restricted to $H^1_r(\mathbb{R}^N)$ are true critical points in $H^1(\mathbb{R}^N)$. Therefore, from now on, we will directly define $J$ on $H^1_r(\mathbb{R}^N)$.

Before stating our assumptions on $M$ and the main results of this paper, we would like to mention the closely related works of Azzollini, d'Avenia and Pomponio \cite{Az11} and Lu \cite{Lu15-1}. To the best of our knowledge, it seems that only articles \cite{Az11} and \cite{Lu15-1} have considered the multiplicity of solutions for such problem under the very general assumptions on $f$.

In \cite{Az11}, under the same very general assumptions on $f$ as above, Azzollini, d'Avenia and Pomponio considered a suitable perturbation of $I$, namely
\begin{equation*}
I_q(u):=I(u)+qR(u),~~~~u\in H^1(\mathbb{R}^N),
\end{equation*}
where $I$ is given by \eqref{equ1.3}, $q>0$ is a positive parameter, $R:H^1(\mathbb{R}^N)\to\mathbb{R}$ and $N\geq3$. The authors supposed that $R=\Sigma^k_{i=1}R_i$ and, for each $i=1,\cdots,k$, the functional $R_i$ satisfies certain suitable assumptions and the following condition:
\begin{itemize}
  \item[$(R2)$] There exists $\delta_i>0$ such that, for any $u\in H^1(\mathbb{R}^N)$, we have
                \begin{equation*}
                      R'_i(u)[u]\leq C \|u\|^{\delta_i}_{H^1(\mathbb{R}^N)}.
                \end{equation*}
\end{itemize}
By a suitable combination of the method described in \cite{Hi10} and a certain truncation argument, they established an abstract theorem which claims the existence of (at least) $n$ distinct critical points of $I_q$ for every $n\in \mathbb{N}$ and $q\in (0,q_n)$, where $q_n>0$ is a suitable positive constant depending on $n$. As an application, in the case where $N\geq3$ and $M(t)=a+bt$ with $a,b>0$, they treated Problem $(\mathcal{KT})$ and obtained finitely many distinct radial solutions for sufficiently small $b>0$. For another application to the nonlinear Schrodinger-Maxwell system, we refer the reader to \cite{Az11}.

We note that the truncation argument explored in \cite{Az11} is important to the proof of the abstract existence result. Actually, the truncation argument is not only used to construct a suitable modified functional of $I_q$, which satisfies the symmetric mountain pass geometry, but also, together with the method described in \cite{Hi10},  plays a vital role in obtaining (at least) $n$ distinct particular Palais-Smale sequences which are bounded for every $n\in \mathbb{N}$ and $q\in (0,q_n)$. Thus, it is interesting to ask the question whether, at least for Problem $(\mathcal{KT})$ in the case where $M(t)=a+bt$ with $a,b>0$ and $N\geq3$, it is possible to prove the multiple result by some suitable arguments, e.g. variational methods, but without using a truncation technique as in \cite{Az11}.

In the more recent paper \cite{Lu15-1}, by means of a scaling argument based on an idea of Azzollini \cite{Az12,Az15} and a new description of the critical values, we investigated the following Kirchhoff Problem
\begin{equation}\tag{$\mathcal{K}$}
\left\{
\begin{aligned}
&-\left(a+b\int_{\mathbb{R}^N}|\nabla{u}|^2\right)\Delta{u}= f(u)~~\text{in}~\mathbb{R}^N,\\
&u\in H^1(\mathbb{R}^N),~~~~u\not\equiv0~~\text{in}~\mathbb{R}^N,\\
\end{aligned}
\right.
\end{equation}
where $a\geq0$, $b>0$ and $N\geq1$. When $N\geq2$, under some suitable conditions on the values of the nonnegative parameters $a$ and $b$ if necessary and the assumptions $(f_0),(f_2),(f_3)$ and $(f'_1)$ on $f$, certain multiplicity results for $(\mathcal{K})$ were obtained as partial results in that paper. In particular, we obtained infinitely many distinct radial solutions in \cite{Lu15-1} for any $a\geq0$ and $b>0$ fixed when $N=2,3$. We note here that \cite{Lu15-1} not only answers the question we raised above in the affirmative from the non-variational point of view, but also extends the result of Azzollini, d'Avenia and Pomponio in \cite{Az11} concerning the existence of multiple solutions to $(\mathcal{K})$.

As pointed out in \cite{Lu15-1}, it is natural to know whether, at least for the non-degenerate case $a>0$, one can still obtain the multiplicity results for $(\mathcal{K})$ via variational methods. So far, this question has a positive answer for the case $N\geq4$ by the early work \cite{Az11} of Azzollini, d'Avenia and Pomponio. However, this question is still open for the cases $N=2,3$, where, in fact, Problem $(\mathcal{K})$ possesses infinitely many distinct radial solutions that .

Motivated by the articles \cite{Az11,Hi10,Lu15-1} and the questions we raised above, by making some suitable assumptions on $M$, we shall show the existence of infinitely many distinct radial solutions for Problem $(\mathcal{KT})$ as our first result of this paper. For this purpose, we make the hypotheses on the function $M$ as follows:
\begin{itemize}
  \item[$(M_1)$] There exists $m_0>0$ such that $M(t)\geq m_0$ for any $t\geq0$.
  \item[$(M_2)$] Let $\widehat{M}(t):=\int^t_0M(\tau)d\tau$. Then we have
  \begin{equation*}
  \underset{t\to +\infty}{\liminf}~\left[\widehat{M}(t)-\left(1-\frac{2}{N}\right)M(t)t\right]=+\infty.
  \end{equation*}
  \item[$(M_3)$] We assume $$\underset{t\to +\infty}{\lim}\frac{M(t)}{t^{\frac{2}{N-2}}}=0.$$
\end{itemize}

Now, our first result of the present paper can be stated as follows.
\begin{theorem}
Assume $N\geq2$ and that $f$ satisfies $(f_0),(f_2),(f_3)$ and $(f'_1)$. Suppose $(M_1)$ when $N=2$ and $(M_1)-(M_3)$ when $N\geq3$. Then Problem $(\mathcal{KT})$ has infinitely many distinct (possibly sign-changing) radially symmetric solutions.
\end{theorem}
\begin{remark}
All the solutions that we obtain in Theorem 1.1 are characterized by the symmetric mountain pass minimax argument in $H^1_r(\mathbb{R}^N)$.
\end{remark}
\begin{remark}
In our proof of Theorem 1.1, a truncation argument similar to that in \cite{Az11} would and should be avoided; since, if not, in general it seems to be difficult or even impossible to get infinitely many distinct solutions. This can be seen as another reason why we try to find solutions of Problem $(\mathcal{KT})$ through the non-modified functional $J$ directly.
\end{remark}
\begin{remark}
When $N=2$, under the same assumptions of Theorem 1.1, Figueiredo, Ikoma and J\'{u}nior have obtained a least energy solution of $(\mathcal{KT})$ in the early work \cite{Fi14}. In that paper, under certain suitable conditions on $M$ which are stronger than $(M_1)-(M_3)$, the existence result of a least energy solutions to $(\mathcal{KT})$ was also established for $N\geq3$. Our Theorem 1.1 here can be viewed as a natural extension of \cite{Fi14}.
\end{remark}

Next, when $N\geq3$, for a suitable class of non-degenerate functions $M$ which may not satisfy the hypothesis $(M_3)$, we establish the following weaker multiplicity result, which claims the existence of finitely many radial solutions to $(\mathcal{KT})$.

\begin{theorem}
Assume that $M(t)=m_0+q \lambda(t)$ with $q>0$ and $\lambda\in C(\mathbb{R}_+,\mathbb{R}_+)$, $N\geq3$ and that $f$ satisfies $(f_0),(f_2),(f_3)$ and $(f'_1)$. Besides, suppose that either $(M_2)$ or the following is satisfied:
\begin{itemize}
  \item[$(M'_2)$] We have the inequality
  \begin{equation*}
  \underset{t\to +\infty}{\limsup}~\left[\widehat{M}(t)-\left(1-\frac{2}{N}\right)M(t)t\right]\leq0,
  \end{equation*}
\end{itemize}
Then, for any $n\in\mathbb{N}$, there exists a positive constant $q_n>0$ such that Problem $(\mathcal{KT})$ has at least $n$ distinct (possibly sign-changing) radial solutions for any $q\in(0,q_n)$.
\end{theorem}
\begin{remark}
All the solutions that we obtain in Theorem 1.2 are characterized by the symmetric mountain pass minimax argument in $H^1_r(\mathbb{R}^N)$.
\end{remark}
\begin{remark}
The fact that, under the assumptions of Theorem 1.2, we obtain only finitely many nontrivial solutions for sufficiently small $q>0$ is not surprising and we can hardly expect more. Actually, the function $M(t)=m_0+qt^{\frac{2}{N-2}}$ with $m_0,q>0$ satisfies $(M_2)$. However, in this case, Theorem A.1 in the paper \cite{Fi14} by Figueiredo, Ikoma and J\'{u}nior showed the nonexistence of nontrivial solution for large enough $q>0$. In addition, by repeating certain arguments explored in \cite{Lu15-1} for the proof of Theorem 1.2, Item $(ii)$ in that paper, we can only show that more and more distinct solutions of $(\mathcal{KT})$ exist as $q\to0^+$. It seems to be difficult or even impossible to get infinitely many distinct solutions of $(\mathcal{KT})$ for sufficiently small but fixed $q>0$. Thus, the conclusion of Theorem 1.2 seems to be the best possible result we could hope for when $M$ does indeed not satisfy the hypothesis $(M_3)$.
\end{remark}
\begin{remark}
It is not difficult to see that the hypothesis $(R2)$ is not always satisfied for functions $M$ which verify the assumptions of Theorem 1.2. For example, let $m_0=1$ and $\lambda(t)=\frac{1}{2}(e^t-1)$, that is $M(t)=1+\frac{q}{2}(e^t-1)$, then a straightforward computation shows that such $M$ satisfies assumption $(M'_2)$. However, in this case, $(R2)$ is not satisfied due to the fact that, for any $\delta>0$ and $u\in H^1(\mathbb{R}^N)\setminus\{0\}$, we have
\begin{equation*}
\underset{t\to+\infty}{\lim}\frac{R'(tu)[tu]}{\|tu\|^{\delta}}=\frac{\|\nabla u\|^2_2}{2\|u\|^{\delta}}\underset{t\to+\infty}{\lim}\left(e^{t^2\|\nabla u\|^2_2}-1\right)t^{2-\delta}=+\infty.
\end{equation*}
Thus, our Theorem 1.2 can not be obtained by applying the abstract result given by \cite{Az11} directly and the arguments there are also not valid here.
\end{remark}

As a consequence of Theorems 1.1 and 1.2, we have the following result:
\begin{corollary}
Assume $a>0$ fixed, $b>0$, $N\geq2$ and that $f$ satisfies $(f_0),(f_2),(f_3)$ and $(f'_1)$. Then the following statements hold.
\begin{itemize}
  \item[~~$(i)$] If $N=2,3$, Problem $(\mathcal{K})$ has infinitely many distinct radially symmetric solutions for any $b>0$, which are characterized by the symmetric mountain pass minimax argument in $H^1_r(\mathbb{R}^N)$.
  \item[~$(ii)$] If $N\geq4$, for any $n\in\mathbb{N}$, there exists a positive constant $b_n>0$ such that Problem $(\mathcal{K})$ has at least $n$ distinct radially symmetric solutions for any $b\in(0,b_n)$. Moreover, all the solutions are characterized by the symmetric mountain pass minimax argument in $H^1_r(\mathbb{R}^N)$.
\end{itemize}
\end{corollary}
\begin{remark}
As we can see in Sections 3 and 4, the proofs of Theorems 1.1 and 1.2 are all based on a certain variational method described in \cite{Hi10} but without a truncation argument similar to that in \cite{Az11}. Since Corollary 1.1 follows directly from Theorems 1.1 and 1.2, we thus answer the first question we raised above in the affirmative again from the variational point of view and address the second problem we raised above in the remaining cases $N=2,3$. As a by-product, in the case $a>0$ and $N\geq4$, we provide another variational proof of the multiple result of $(\mathcal{K})$ through the non-modified functional $J$, which is different from that in \cite{Az11}.
\end{remark}

The remaining part of this paper is organized as follows. In Section 2, an auxiliary problem is constructed in the spirit of \cite{Hi10} and the corresponding conclusions are shown at the same time. With the aid of the method described in \cite{Hi10} and the conclusions in Section 2, the proofs of Theorems 1.1 and 1.2 are completed in Sections 3 and 4 respectively. Lastly, in Section 5, the non-variational proofs are presented which actually provide us a better understanding of the multiplicity results.

\section{The auxiliary problem and its result}
In this section, we shall construct an auxiliary problem in the spirit of \cite{Hi10}, which will play an important role in the proofs of the main results of this paper. To be more precise, it will be proved that there is a sequence of positive critical values $\{e_n\}^{+\infty}_{n=1}$ corresponding to the auxiliary problem which is divergent to infinity. This fact allows us to prove the multiplicity results for our original problem $(\mathcal{KT})$ based on the level sets argument.

Following \cite{Hi10}, we set
\begin{equation*}
\omega:=-\frac{1}{2}\underset{t\to0}{\limsup}\frac{f(t)}{t}\in(0,+\infty)
\end{equation*}
and equip $H^1_r(\mathbb{R}^N)$ with the norm $\|\cdot\|:=\left(m_0\|\nabla \cdot\|^2_2+\omega\|\cdot\|^2_2\right)^{\frac{1}{2}}$.

Consider $p_0\in \left(1,\frac{N+2}{N-2}\right)$ if $N\geq3$, $p_0\in (1,+\infty)$ if $N=2$ and set
\begin{equation*}
\begin{split}
&h(t):=\left\{
\begin{aligned}
&\text{max}\{\omega t+f(t),0\},~&\text{for}~t\geq0,\\
&-h(-t),&\text{for}~t<0,\\
\end{aligned}
\right.
~~~~~~H(t):=\int^t_0h(\tau)d\tau,\\
&\overline{h}(t):=\left\{
\begin{aligned}
&t^{p_0}\underset{0<\tau\leq t}{\max}\frac{h(\tau)}{\tau^{p_0}},~~~~~~~&\text{for}~t>0,\\
&0,&\text{for}~t=0,\\
&-\overline{h}(-t),&\text{for}~t<0,\\
\end{aligned}
\right.
~~~~~~\overline{H}(t):=\int^t_0\overline{h}(\tau)d\tau.
\end{split}
\end{equation*}
Then, the functions $h,\overline{h},H$ and $\overline{H}$ satisfy the properties stated in Lemmas 2.1-2.3 below.
\begin{lemma}
The following statements hold:
\begin{itemize}
  \item[$~~(i)$] For all $t\geq0$, we have $\omega t+f(t)\leq h(t)\leq \overline{h}(t)$.
  \item[$~(ii)$] For all $t\geq0$, we have $\overline{h}(t),h(t)\geq0$.
  \item[$(iii)$] There exists $\delta>0$ such that $h(t)=\overline{h}(t)=0$ for all $t\in[0,\delta]$
  \item[$(iv)$] There exists $\xi>0$ such that $0<h(\xi)\leq \overline{h}(\xi)$.
  \item[~$(v)$] The map $t\mapsto \frac{\overline{h}(t)}{t^{p_0}}$ is non-decreasing in $t\in (0,+\infty)$.
  \item[$(vi)$] The functions $h$ and $\overline{h}$ satisfy $(f_2)$.
\end{itemize}
\end{lemma}

\begin{lemma}
The following statements hold:
\begin{itemize}
  \item[$~~(i)$] For all $t\geq0$, we have $\frac{1}{2}\omega t^2+F(t)\leq H(t)\leq \overline{H}(t)$.
  \item[$~(ii)$] For all $t\geq0$, we have $\overline{H}(t),H(t)\geq0$.
  \item[$(iii)$] There exists $\delta>0$ such that $H(t)=\overline{H}(t)=0$ for all $t\in[0,\delta]$
  \item[$(iv)$] It holds that $\overline{H}(\zeta)-\frac{1}{2}\omega \zeta^2>0$.
  \item[~$(v)$] For all $t\in \mathbb{R}$, $0\leq (p_0+1)\overline{H}(t)\leq t\overline{h}(t)$.
  \item[$(vi)$] The functions $H$ and $\overline{H}$ satisfy
  \begin{equation*}
  \begin{split}
  &\underset{|t|\to+\infty}{\lim}\frac{H(t)}{t^{\frac{2N}{N-2}}}=\underset{|t|\to+\infty}{\lim}\frac{\overline{H}(t)}{t^{\frac{2N}{N-2}}}=0,~~~~~~\text{when}~N\geq3,\\
  &\underset{|t|\to+\infty}{\lim}\frac{H(t)}{e^{\alpha t^2}}=\underset{|t|\to+\infty}{\lim}\frac{\overline{H}(t)}{e^{\alpha t^2}}=0,~~~~~~\text{for any}~\alpha>0~\text{when}~N=2.
  \end{split}
  \end{equation*}
\end{itemize}
\end{lemma}

\begin{lemma}
Let $N\geq2$ and suppose that $\{u_j\}^{+\infty}_{j=1}\subset H^1_r(\mathbb{R}^N)$ converges to $u\in H^1_r(\mathbb{R}^N)$ weakly in $H^1_r(\mathbb{R}^N)$. Then
\begin{itemize}
\item[~~$(i)$] $\int_{\mathbb{R}^N}H(u_j)\to \int_{\mathbb{R}^N}H(u)$ and $\int_{\mathbb{R}^N}\overline{H}(u_j)\to \int_{\mathbb{R}^N}\overline{H}(u)$.
\item[~$(ii)$] $h(u_j)\to h(u)$ and $\overline{h}(u_j)\to \overline{h}(u)$ strongly in $(H^1_r(\mathbb{R}^N))^{-1}$.
\end{itemize}
\end{lemma}

Now, we can construct the auxiliary problem as follow:
\begin{equation}\tag{$\mathcal{A}$}
\left\{
\begin{aligned}
&-m_0\Delta{u}+\omega u= \overline{h}(u),~~\text{in}~\mathbb{R}^N,\\
&u\in H^1_r(\mathbb{R}^N),~~~~u\not\equiv0~~~\text{in}~\mathbb{R}^N,\\
\end{aligned}
\right.
\end{equation}
where $N\geq2$, $m_0>0$ given by $(M_1)$, $\omega>0$ and $\overline{h}\in C(\mathbb{R},\mathbb{R})$ defined as above. It is not difficult to see that the corresponding functional of $(\mathcal{A})$ given by
\begin{equation*}
K(u):=\frac{1}{2}\left\|u\right\|^2-\int_{\mathbb{R}^N}\overline{H}(u)
\end{equation*}
is well-defined on $H^1_r(\mathbb{R}^N)$ and of class $C^1$. Moreover, as stated in the next lemma, $K$ has the geometry of the Symmetric Mountain Pass theorem and satisfies the Palais-Smale compactness condition. In what follows, we set $\mathbb{D}_n:=\left\{\sigma=(\sigma_1,\cdots,\sigma_n)\in \mathbb{R}^n~|~|\sigma|\leq1\right\}$ and $\mathbb{S}^{n-1}:=\partial\mathbb{D}_n.$

\begin{lemma}
The functional $K$ satisfies the following properties.
\begin{itemize}
  \item[~~$(i)$] There exist $r>0$ and $\rho>0$ such that
                 \begin{equation*}
                 \begin{split}
                        K(u)\geq0~~~~\text{for any}~u\in H^1_r(\mathbb{R}^N)~\text{with}~\|u\|\leq r,\\
                        K(u)\geq\rho~~~~\text{for any}~u\in H^1_r(\mathbb{R}^N)~\text{with}~\|u\|= r.
                 \end{split}
                 \end{equation*}
  \item[~$(ii)$] For every $n\in \mathbb{N}$, there exists an odd continuous mapping $\gamma_{0n}:\mathbb{S}^{n-1}\to H^1_r(\mathbb{R}^N)$ such that
    \begin{equation*}
    K\left(\gamma_{0n}(\sigma)\right)<0~~~~\text{for all}~\sigma\in \mathbb{S}^{n-1}.
    \end{equation*}
    \item[$(iii)$] The Palais-Smale compactness condition holds.
\end{itemize}
\end{lemma}

Due to Item $(ii)$ of Lemma 2.4, for every $n\in \mathbb{N}$, we can define a family of mappings $\Gamma_n$ by
\begin{equation}\label{equ2.1}
  \Gamma_n:=\left\{\gamma\in C(\mathbb{D}_n,H^1_r(\mathbb{R}^N))~|~\gamma~\text{is odd and}~\gamma(\sigma)=\gamma_{0n}(\sigma)~\text{on}~\sigma\in \mathbb{S}^{n-1}\right\},
\end{equation}
which is nonempty since
\begin{equation*}
\gamma_n(\sigma):=\left\{
\begin{aligned}
&|\sigma|\gamma_{0n}\left(\frac{\sigma}{|\sigma|}\right),~~~~\text{for}~\sigma\in \mathbb{D}_n\setminus\{0\},\\
&0,~~~~~~~~~~~~~~~~~~~~\text{for}~\sigma=0,\\
\end{aligned}
\right.
\end{equation*}
belongs to $\Gamma_n$. Thus, the symmetric mountain pass values of $K$ defined by
 \begin{equation*}
e_n:=\underset{\gamma\in \Gamma_n}{\inf}\underset{\sigma\in \mathbb{D}_n}{\max}~K(\gamma(\sigma))
 \end{equation*}
for any $n\in \mathbb{N}$, are all meaningful. Moreover, we have the following theorem.
\begin{theorem}
The following statements hold.
\begin{itemize}
  \item[$~(i)$] For every $n\in \mathbb{N}$, $e_n$ is a critical value of $K$ and $e_n\geq\rho>0$.
  \item[$(ii)$] $e_n\to +\infty$ as $n\to+\infty$.
\end{itemize}
\end{theorem}
\begin{remark}
All of the conclusions stated in this section and their proofs can be found in \cite{Hi10}. For ease of exposition and completeness of this paper, it is better to outline the necessary conclusions that we need.
\end{remark}

\section{Proof of Theorem 1.1}
In this section, we shall give the detailed proof of Theorem 1.1. Before going further, we would like to point out that assumption $(M_3)$ is almost necessary when it comes to obtaining infinitely many distinct solutions to $(\mathcal{KT})$ in the case $N\geq3$, see Remark 1.5 in Section 1. On the other hand, as we can see below, actually assumption $(M_3)$ is important to verifying the symmetric mountain pass geometry of $J$ for every $n\in \mathbb{N}$ and, together with assumptions $(M_1)$ and $(M_2)$, is also sufficient to establish the existence result of infinite many distinct solutions to $(\mathcal{KT})$.

\subsection{Symmetric mountain pass geometry of $J$}

\begin{lemma}
The functional $J$ also satisfies Items $(i)$ and $(ii)$ of Lemma 2.4 in Section 2.
\end{lemma}
\proof  In terms of Item $(i)$ of Lemma 2.2 and $(M_1)$, we have
\begin{equation}\label{equ3.1}
J(u)\geq K(u)~~~~ \text{for all}~ u\in H^1_r(\mathbb{R}^N),
\end{equation}
which implies that Item $(i)$ of Lemma 2.4 is applied to $J$.

For every $n\in \mathbb{N}$, arguing as in Theorem 10 of \cite{Be83-2}, an odd and continuous map $\pi_n:\mathbb{S}^{n-1}\to H^1_r(\mathbb{R}^N)$ is defined such that
\begin{equation*}
0\notin \pi_n\left(\mathbb{S}^{n-1}\right)~~~~\text{and}~~~~\int_{\mathbb{R}^N}F\left(\pi_n(\sigma)\right)\geq1,~~\text{for all}~\sigma\in\mathbb{S}^{n-1}.
\end{equation*}
It is easy to see that, for every $n\in\mathbb{N}$, there exists $\alpha_n>0$ such that
\begin{equation*}
\left\|\nabla \pi_n(\sigma)\right\|^2_2\leq \alpha_n,~~\text{for all}~\sigma\in \mathbb{S}^{n-1}.
\end{equation*}
For every $n\in\mathbb{N}$ and any $\sigma\in\mathbb{S}^{n-1}$, setting $\beta^t_n(\sigma)(x):=\pi_n(\sigma)(t^{-1}x)$, we have
\begin{equation*}
\begin{split}
J\left(\beta^t_n(\sigma)\right)
&=\frac{1}{2}\widehat{M}\left(t^{N-2}\|\nabla \pi_n(\sigma)\|^2_2\right)-t^N\int_{\mathbb{R}^N}F\left(\pi_n(\sigma)\right)\\
&\leq\frac{1}{2}\widehat{M}\left(t^{N-2}\alpha_n\right)-t^N=:g_n(t).
\end{split}
\end{equation*}
 When $N=2$, it is clear that $g_n(t_n)<0$ for sufficiently large $t_n>0$. When $N\geq3$, in terms of $(M_3)$, there also exists a sufficiently large $t_n>0$ such that
\begin{equation*}
g_n(t_n)=t^N_n\left(\frac{1}{2}\frac{\widehat{M}(s_n)}{s^{\frac{N}{N-2}}_n}\alpha^{\frac{N}{N-2}}_n-1\right)<0,
\end{equation*}
where $s_n:=t^{N-2}_n\alpha_n>0$. Thus the proof is completed by redefining $\gamma_{0n}:=\beta^{t_n}_n$.~~$\square$

Now, for $n\in \mathbb{N}$, we can defined the symmetric mountain pass values $d_n$ of $J$:
\begin{equation*}
d_n:=\underset{\gamma\in \Gamma_n}{\inf}\underset{\sigma\in \mathbb{D}_n}{\max}~J(\gamma(\sigma)),
\end{equation*}
where $\Gamma_n$ is given by \eqref{equ2.1}. In view of \eqref{equ3.1} and Theorem 2.1, we have that
\begin{equation*}
d_n\geq e_n\geq \rho>0~~~~\text{and}~~~~d_n\to+\infty~~\text{as}~n\to+\infty.
\end{equation*}
It is easy to see that the proof of Theorem 1.1 is completed if we can prove that, for every $n\in\mathbb{N}$, the value $d_n$ defined above is a critical value of $J$.

For every $n\in\mathbb{N}$, by Ekeland's principle, we can find a Palais-Smale sequence $\{u_j\}^{+\infty}_{j=1}$ at level $d_n$, that is, $\{u_j\}^{+\infty}_{j=1}$ satisfies
\begin{equation}\label{equ3.2}
J(u_j)\to d_n~~~~\text{and}~~~~J'(u_j)\to 0~\text{in}~(H^1_r(\mathbb{R}^N))^{-1},~~~~\text{as}~j\to+\infty.
\end{equation}
However, merely under the condition \eqref{equ3.2}, it seems difficult to show the existence of a strongly convergent subsequence and even the boundedness of $\{u_j\}^{+\infty}_{j=1}$ in $H^1_r(\mathbb{R}^N)$. Inspired by \cite{Hi10}, by introducing an auxiliary functional, we find a Palais-Smale sequence that "almost" satisfies the Poho\v{z}aev identity associated to $(\mathcal{KT})$, which makes it possible for us to overcome these difficulties.

In the following subsection, based on the key idea above, we will show that $d_n$ is indeed a critical value of $J$ for every $n\in\mathbb{N}$.

\subsection{Auxiliary functional $\Phi(\theta,u)$ and conclusion}

Analogously to \cite{Hi10}, we equip a standard product norm $\|(\theta,u)\|_{\mathbb{R}\times H^1_r}:=\left(\theta^2+\|u\|^2\right)^{\frac{1}{2}}$
to the augmented space $\mathbb{R}\times H^1_r(\mathbb{R}^N)$ and define the auxiliary functional
\begin{equation*}
\Phi(\theta,u):=\frac{1}{2}\widehat{M}\left(e^{(N-2)\theta}\int_{\mathbb{R}^N}|\nabla u|^2\right)-e^{N\theta}\int_{\mathbb{R}^N}F(u).
\end{equation*}
It is easy to conclude that $\Phi$ is of class $C^1$ and
\begin{equation}\label{equ3.3}
\Phi(\theta,u(x))=J\left(u\left(e^{-\theta}x\right)\right)~~~~\text{for all}~ \theta\in\mathbb{R}~ \text{and} ~u\in H^1_r(\mathbb{R}^N).
\end{equation}
In particular, $\Phi(0,u)=J\left(u\right)$ for all $u\in H^1_r(\mathbb{R}^N)$. We denote its derivative as $\Phi':=\left(\partial_\theta\Phi,\partial_u\Phi\right)$ with
\begin{equation*}
\partial_\theta\Phi(\theta,u)=\frac{N-2}{2}M\left(e^{(N-2)\theta}\int_{\mathbb{R}^N}|\nabla u|^2\right)e^{(N-2)\theta}\int_{\mathbb{R}^N}|\nabla u|^2-Ne^{N\theta}\int_{\mathbb{R}^N}F(u)
\end{equation*}
and
\begin{equation*}
\partial_u\Phi(\theta,u)[v]=M\left(e^{(N-2)\theta}\int_{\mathbb{R}^N}|\nabla u|^2\right)e^{(N-2)\theta}\int_{\mathbb{R}^N}\nabla u\cdot\nabla v-e^{N\theta}\int_{\mathbb{R}^N}f(u)v,
\end{equation*}
for all $v\in H^1_r(\mathbb{R}^N)$.

For every $n\in \mathbb{N}$, we define the class
\begin{equation*}
\overline{\Gamma}_n:=\left\{
\begin{aligned}
\overline{\gamma}\in C(\mathbb{D}_n,\mathbb{R}\times H^1_r(\mathbb{R}^N))~\left|~
\begin{aligned}
&\overline{\gamma}(\sigma)=\left(\theta(\sigma),\eta(\sigma)\right)~
\text{satisfies}~\\
&\left(\theta(-\sigma),\eta(-\sigma)\right)=\left(\theta(\sigma),-\eta(\sigma)\right),&\forall\sigma\in \mathbb{D}_n,&\\
&\left(\theta(\sigma),\eta(\sigma)\right)=\left(0,\gamma_{0n}(\sigma)\right),&\forall\sigma\in \mathbb{S}^{n-1}.&
\end{aligned}
\right.
\end{aligned}
\right\},
\end{equation*}
where $\gamma_{0n}$ is given in Item $(ii)$ of Lemma 2.4. In terms of the nonemptyness of $\Gamma_n$ and the fact that $\left\{(0,\gamma)~|~\gamma\in \Gamma_n\right\}\subset\overline{\Gamma}_n$, we conclude that $\overline{\Gamma}_n$ is nonempty, the minimax value $\overline{d}_n$ of $\Phi$ given by
\begin{equation*}
\overline{d}_n:=\underset{\overline{\gamma}\in \overline{\Gamma}_n}{\inf}\underset{\sigma\in D_n}{\max}\Phi\left(\overline{\gamma}\left(\sigma\right)\right)
\end{equation*}
is well-defined and $\overline{d}_n\leq d_n$. On the other hand, for any given $\overline{\gamma}(\sigma)=\left(\theta(\sigma),\eta(\sigma)\right)\in\overline{\Gamma}_n$, setting $\gamma(\sigma)(x)=\eta(\sigma)\left(e^{-\theta(\sigma)}x\right)$, we can verify that $\gamma(\sigma)\in\Gamma_n$ and, by \eqref{equ3.3}, $I(\gamma(\sigma))=\Phi(\overline{\gamma}(\sigma))$ for any $\sigma\in \mathbb{D}_n$, which imply that $\overline{d}_n\geq d_n$. Thus, we conclude the following lemma.
\begin{lemma}
For all $n\in \mathbb{N}$, we have $\overline{d}_n=d_n$.
\end{lemma}

Based on Lemma 3.2, arguing as the proof of Proposition 4.2 in \cite{Hi10}, we have the following lemma:
\begin{lemma}
For every $n\in \mathbb{N}$, there exists a sequence $\left\{\left(\theta_j,u_j\right)\right\}^{+\infty}_{j=1}\subset \mathbb{R}\times H^1_r(\mathbb{R}^N)$ such that
\begin{itemize}
\item[~~$(i)$] $\theta_j\to0$,
\item[~$(ii)$] $\Phi(\theta_j,u_j)\to d_n$,
\item[$(iii)$] $\partial_u\Phi\left(\theta_j,u_j\right)\to0$ strongly in $(H^1_r(\mathbb{R}^N))^{-1}$,
\item[$(iv)$] $\partial_\theta\Phi \left(\theta_j,u_j\right)\to 0$.
\end{itemize}
\end{lemma}

\begin{lemma}
Let $\left\{\left(\theta_j,u_j\right)\right\}^{+\infty}_{j=1}$ be the sequence given by Lemma 3.3. Then $\{u_j\}^{+\infty}_{j=1}$ is bounded and has a strongly convergent subsequence in $H^1_r(\mathbb{R}^N)$.
\end{lemma}
\proof  For the sake of clarity, we divide the proof into three claims. We shall prove the boundedness of $\{\|\nabla u_j\|^2_2\}^{+\infty}_{j=1}$ in Claim 1, complete the boundedness of $\{u_j\}^{+\infty}_{j=1}$ in $H^1_r(\mathbb{R}^N)$ by showing, in Claim 2, that $\{\|u_j\|^2_2\}^{+\infty}_{j=1}$ is bounded and conclude the existence of a strongly convergent subsequence in Claim 3.

\bigskip
\noindent
\textbf{Claim 1.} {\it The sequence $\{\|\nabla u_j\|^2_2\}^{+\infty}_{j=1}$ is bounded.}

In view of Items $(ii)$ and $(iv)$ of Lemma 3.3, setting $\mu_j:=e^{(N-2)\theta_j}\int_{\mathbb{R}^N}|\nabla u_j|^2$, we have
\begin{equation*}
\widehat{M}\left(\mu_j\right)-\left(1-\frac{2}{N}\right)M\left(\mu_j\right)\mu_j
=2 \Phi(\theta_j,u_j)-\frac{2}{N}\partial_\theta\Phi(\theta_j,u_j)
=2d_n+o_j(1).
\end{equation*}
Thus, in association with the Item $(i)$ of Lemma 3.3, we conclude the boundedness of $\{\|\nabla u_j\|^2_2\}^{+\infty}_{j=1}$ from $(M_1)$ when $N=2$ and $(M_2)$ when $N\geq3$ respectively.

\bigskip
\noindent
\textbf{Claim 2.} {\it The sequence $\{\|u_j\|^2_2\}^{+\infty}_{j=1}$ is bounded and then, by Claim 1, $\{u_j\}^{+\infty}_{j=1}$ is bounded in $H^1_r(\mathbb{R}^N)$.}

Arguing by contradiction, let us assume that, up to a subsequence, $\|u_j\|_2\to+\infty$. For every $j\in\mathbb{N}$, set $t_j:=\|u_j\|^{-\frac{2}{N}}_2$ and $v_j(\cdot):=u_j(t^{-1}_j\cdot)$. Then,
\begin{equation}\label{equ3.4}
t_j\to 0, ~~~~\text{as}~ j\to+\infty.
\end{equation}
By some simple calculations, we have
\begin{equation}\label{equ3.5}
\nabla v_j(\cdot)=t^{-1}_j\nabla u_j(t^{-1}_j\cdot),~~~~\|\nabla v_j \|^2_2=t^{N-2}_j\|\nabla u_j\|^2_2,~~~~\|v_j\|^2_2=1,
\end{equation}
which imply the boundedness of $\{v_j\}$ in $H^1_r(\mathbb{R}^N)$ with the aid of Claim 1 and \eqref{equ3.4}. Without loss of generality, up to a subsequence, we may assume that $v_j\rightharpoonup v_0$ in $H^1_r(\mathbb{R}^N)$. Set $\varepsilon_j:=\|\partial_u\Phi(\theta_j,u_j)\|_{(H^1_r(\mathbb{R}^N))^{-1}}$, with the help of \eqref{equ3.5} and Item $(i)$ of Lemma 2.1, some calculations show that
\begin{equation*}
\begin{split}
\omega e^{N\theta_j}\leq&~ M\left(e^{(N-2)\theta_j}\int_{\mathbb{R}^N}|\nabla u_j|^2\right)e^{(N-2)\theta_j}t^2_j\int_{\mathbb{R}^N}|\nabla v_j|^2+\omega e^{N\theta_j}\int_{\mathbb{R}^N}v^2_j\\
=&~\partial_u\Phi(\theta_j,u_j)\left[t^N_ju_j\right]+e^{N\theta_j}\int_{\mathbb{R}^N}\left(f(v_j)+\omega v_j\right)v_j\\
\leq&~ \varepsilon_j\left(m_0t^{N}_j\|\nabla u_j\|^2_2+\omega\right)^{\frac{1}{2}}+e^{N\theta_j}\int_{\mathbb{R}^N}h(v_j)v_j.
\end{split}
\end{equation*}
Then, by \eqref{equ3.4}, Claim 1, Item $(ii)$ of Lemma 2.3 and Items $(i)$ and $(iii)$ of Lemma 3.3, we have
\begin{equation*}
0<\omega\leq\int_{\mathbb{R}^N}h(v_0)v_0,
\end{equation*}
which implies $v_0\not\equiv0$.

On the other hand, let $\varphi\in H^1_r(\mathbb{R}^N)$ be a function with compact support and, for every $j\in\mathbb{N}$, set $\psi_j(\cdot):=\varphi(t_j\cdot)$. With the aid of the fact that $v_j\rightharpoonup v_0$ in $H^1_r(\mathbb{R}^N)$, Items $(i)$ and $(iii)$ of Lemma 3.3, Claim 1 and \eqref{equ3.4}, we have
\begin{equation*}
\begin{aligned}
\left|\int_{\mathbb{R}^N}f(v_0)\varphi\right|=&~\left|e^{N\theta_j}\int_{\mathbb{R}^N}f(v_j)\varphi\right|+o_j(1)\\
\leq&~\left|\partial_u\Phi(\theta_j,u_j)\left[t^N_j\psi_j\right]\right|\\
&~~~~~~~~+ M\left(e^{(N-2)\theta_j}\int_{\mathbb{R}^N}|\nabla u_j|^2\right)e^{(N-2)\theta_j}t^2_j\left|\int_{\mathbb{R}^N}\nabla v_j\nabla \varphi\right|+o_j(1)\\
\leq&~ \varepsilon_j\left(m_0t^{2}_j\|\nabla \varphi\|^2_2+\omega\|\varphi\|^2_2\right)^{\frac{1}{2}}+Ct^2_j+o_j(1)\to 0.
\end{aligned}
\end{equation*}
Thus, there holds
\begin{equation*}
\int_{\mathbb{R}^N}f(v_0)\varphi=0,~~~~\text{for any}~ \varphi\in H^1_r(\mathbb{R}^N)~ \text{with compact support},
\end{equation*}
which implies $f(v_0)\equiv0$. However, from $(f'_1)$, it follows that $0$ is an isolated zero point of $f$. In association with the fact that $H^1_r(\mathbb{R}^N)\subset C(\mathbb{R}^N\setminus\{0\})$ and $v_0(x)\to0$ as $|x|\to+\infty$, e.g. see \cite{Be83-1}, we have $v_0\equiv0$, which is a contradiction.

\bigskip
\noindent
\textbf{Claim 3.} {\it The sequence $\{u_j\}^{+\infty}_{j=1}$ has a strongly convergent subsequence in $H^1_r(\mathbb{R}^N)$.}

From Claim 2 and Item $(i)$ of Lemma 3.3, up to a subsequence, we may assume that, when $j$ tends to infinity, $u_j\rightharpoonup u_0$ weakly in $H^1_r(\mathbb{R}^N)$ and
\begin{equation*}
\alpha_j:=M\left(e^{(N-2)\theta_j}\|\nabla u_j\|^2_2\right)\to \alpha_0\in(0,+\infty).
\end{equation*}
Then, by Items $(i)$ and $(iii)$ of Lemma 3.3, it is not difficult to see that $u_0$ satisfies
\begin{equation*}
-\alpha_0\Delta u_0=f(u_0)~~\text{in}~\mathbb{R}^N,
\end{equation*}
which implies
\begin{equation}\label{equ3.6}
\alpha_0\|\nabla u_0\|^2_2=\int_{\mathbb{R}^N}f(u_0)u_0.
\end{equation}

On the other hand, a straightforward computation yields
\begin{equation}\label{equ3.7}
\alpha_je^{(N-2)\theta_j}\|\nabla u_j\|^2_2+\omega e^{N\theta_j}\|u_j\|^2_2=\partial_u \Phi(\theta_j,u_j)[u_j]+ \beta^1_je^{N\theta_j}-\beta^2_je^{N\theta_j}
\end{equation}
where
\begin{equation*}
\beta^1_j:=\int_{\mathbb{R}^N}h(u_j)u_j~~~~\text{and}~~~~\beta^2_j:=\int_{\mathbb{R}^N}\left(h(u_j)u_j-f(u_j)u_j-\omega u^2_j\right).
\end{equation*}
Noting that, by Item $(iii)$ of Lemma 3.3, Item $(ii)$ of Lemma 2.3 and Item $(i)$ of Lemma 2.1 and Fatou's lemma, there hold
\begin{equation}\label{equ3.8}
\underset{j\to+\infty}{\lim}\partial_u \Phi(\theta_j,u_j)[u_j]=0,~~~~\underset{j\to+\infty}{\lim}\beta^1_j=\int_{\mathbb{R}^N}h(u_0)u_0,
\end{equation}
and
\begin{equation}\label{equ3.9}
\underset{j\to+\infty}{\liminf}~\beta^2_j\geq\int_{\mathbb{R}^N}\left(h(u_0)u_0-f(u_0)u_0-\omega u^2_0\right).
\end{equation}
Now, from Item $(i)$ of Lemma 3.3 and \eqref{equ3.6}-\eqref{equ3.9} above, we conclude
\begin{equation*}
\begin{split}
\underset{j\to+\infty}{\limsup}~\left(\alpha_0\|\nabla u_j\|^2_2+\omega\|u_j\|^2_2\right)
&=\underset{j\to+\infty}{\limsup}~\left(\alpha_je^{(N-2)\theta_j}\|\nabla u_j\|^2_2+\omega e^{N\theta_j}\|u_j\|^2_2\right)\\
&\leq\int_{\mathbb{R}^N}\left(f(u_0)u_0+\omega u^2_0\right)\\
&=\alpha_0\|\nabla u_0\|^2_2+\omega\|u_0\|^2_2,
\end{split}
\end{equation*}
which implies that $u_j\to u_0$ in $H^1_r(\mathbb{R}^N)$. Thus the proof of Lemma 3.4 is completed.~~$\square$

\bigskip
\noindent
\textbf{Conclusion}~~ Let $\{(\theta_j,u_j)\}^{+\infty}_{j=1}$ be the sequence given by Lemma 3.3. By Lemma 3.4, we may assume that $u_j\to u_{0n}$ in $H^1_r(\mathbb{R}^N)$. Then, in association with Items $(i)$ and $(iii)$ of Lemma 3.3, it follows that
\begin{equation*}
\Phi(0,u_{0n})=d_n~~~~\text{and}~~~~\partial_u\Phi(0,u_{0n})=0,
\end{equation*}
that is
\begin{equation*}
J(u_{0n})=d_n~~~~\text{and}~~~~J'(u_{0n})=0,
\end{equation*}
Thus, for every $n\in\mathbb{N}$, the symmetric mountain pass value $d_n$ defined in Subsection 3.1 is indeed a critical value of $J$ and, by \eqref{equ3.2}, we complete the proof of Theorem 1.1.~~$\square$

\section{Proof of Theorem 1.2}

In this section, the proof of Theorem 1.2 shall be completed. It is worth pointing out that, due to the loss of assumption $(M_3)$, finding suitable candidate critical values of $J$ becomes the major difficulty that we need to overcome in the proof of Theorem 1.2. Fortunately, as we can see below, we are able to get through this obstacle by Item $(ii)$ of Theorem 2.1 and the non-negativeness of function $\lambda$. As the core of this section, the process of finding suitable candidate critical values will be shown in detail.

For convenience, we rewrite the corresponding functional of $(\mathcal{KT})$ as
\begin{equation*}
\begin{split}
J_q(u):=&~\frac{1}{2}m_0\int_{\mathbb{R}^N}|\nabla u|^2+\frac{q}{2}\Lambda\left(\int_{\mathbb{R}^N}|\nabla u|^2\right)-\int_{\mathbb{R}^N}F(u)\\
=&~I_0(u)+\frac{q}{2}\Lambda\left(\int_{\mathbb{R}^N}|\nabla u|^2\right),
\end{split}
\end{equation*}
where $\Lambda(t):=\int^t_0\lambda(\tau)d\tau$ and $I_0\in C^1(H^1_r(\mathbb{R}^N),\mathbb{R})$ given by
\begin{equation*}
I_0:=\frac{1}{2}m_0\int_{\mathbb{R}^N}|\nabla u|^2-\int_{\mathbb{R}^N}F(u).
\end{equation*}
Apparently, there holds $J_q(u)\geq I_0(u)\geq K(u)$ for all $u\in H^1_r(\mathbb{R}^N)$.

Redefining $\gamma_{0n}$ if necessary, by Lemma 3.1, we have that $I_0$ also satisfies Items $(i)$ and $(ii)$ of Lemma 2.4 in Section 2. It is easy to see that, for such $\gamma_{0n}$, there exist $\alpha_n,\beta_n>0$ such that
\begin{equation*}
I_0(\gamma_{0n}(\sigma))\leq-2\alpha_n<0~~~~\text{and}~~~~\Lambda\left(\int_{\mathbb{R}^N}|\nabla \gamma_{0n}(\sigma)|^2\right)\leq 2\beta_n,~~~~\text{for all}~\sigma\in \mathbb{S}^{n-1}.
\end{equation*}
Let $q\in (0,\alpha_n\beta^{-1}_n]$, then
\begin{equation*}
J_q(\gamma_{0n}(\sigma))=I_0(\gamma_{0n}(\sigma))+\frac{q}{2}\Lambda\left(\int_{\mathbb{R}^N}|\nabla \gamma_{0n}(\sigma)|^2\right)\leq -2\alpha_n+q\beta_n\leq-\alpha_n<0
\end{equation*}
for all $\sigma\in \mathbb{S}^{n-1}$. Therefor, we can define a candidate critical value $c^q_n$ of $J_q$ by
\begin{equation*}
c^q_n:=\underset{\gamma\in \Gamma_n}{\inf}\underset{\sigma\in \mathbb{D}_n}{\max}~J_q(\gamma(\sigma)),
\end{equation*}
where $\Gamma_n$ is given by \eqref{equ2.1}. Obviously, for any $0<q\leq q'\leq \alpha_n\beta^{-1}_n$,
\begin{equation*}
c^{q'}_n\geq c^q_n\geq e_n\geq \rho>0.
\end{equation*}

We claim that, for every $m\in\mathbb{N}$, there exist $\{n_k\}^m_{k=1}\subset\mathbb{N}$ and $q_m>0$ such that, for $q\in(0,q_m]$, the minimax values $\{c^q_{n_k}\}^m_{k=1}$ of $J_q$ given by
\begin{equation*}
c^q_{n_k}:=\underset{\gamma\in \Gamma_{n_k}}{\inf}\underset{\sigma\in \mathbb{D}_{n_k}}{\max}~J_q(\gamma(\sigma)),~~~~k=1,2,\cdots,m
\end{equation*}
are all well-defined and satisfy
\begin{equation*}
0<\rho\leq c^q_{n_1}<\cdots<c^q_{n_k}<\cdots< c^q_{n_m}<+\infty.
\end{equation*}

Actually, let $n_1=1$, $q^{n_1}:=\alpha_{n_1}\beta^{-1}_{n_1}$ and $q\in (0,q^{n_1}]$, it is easy to see that $c^q_{n_1}$ is well defined and
\begin{equation*}
0<\rho\leq e_{n_1}\leq c^q_{n_1}\leq c^{q^{n_1}}_{n_{1}}<+\infty.
\end{equation*}
In view of Item $(ii)$ of Theorem 2.1, there exists $n_2\in\mathbb{N}$ such that $c^{q^{n_1}}_{n_1}<e_{n_2}$. For such $n_2\in\mathbb{N}$, let $q^{n_2}:=\min\{\alpha_{n_2}\beta^{-1}_{n_2},q^{n_1}\}$ and $q\in(0,q^{n_2}]$, we have that $\{c^q_{n_k}\}^2_{k=1}$ are well-defined and satisfy
\begin{equation*}
0<\rho\leq e_{n_1}\leq c^q_{n_1}<e_{n_2}\leq c^q_{n_2}\leq c^{q_{n_2}}_{n_2}<+\infty.
\end{equation*}
Thus, for every fixed $m\in\mathbb{N}$, the desired sequence $\{n_k\}^m_{k=1}\subset\mathbb{N}$ can be obtained by an iterative procedure, and the desired positive number $q_m$ can also be found by letting $q_m:=\underset{1\leq k\leq m}{\min}\{\alpha_{n_k}\beta^{-1}_{n_k}\}$.

It is not difficult to see that, under the assumptions of Theorem 1.2, the arguments explored in Subsection 3.2 are also valid here. This fact means that the candidate critical values of $J_q$ we define above are indeed critical values of $J_q$. Therefor, the proof of Theorem 1.2 is finished.~~$\square$

\section{Non-variational proofs of the multiplicity results}

In this last section, inspired by \cite{Az12,Lu15-1}, we shall present different new proofs of the multiplicity results which are non-variational, simple and fundamental. As we can see below, this gives us natural interpretations of the results we proved in previous sections.

Before going into the details of the non-variational proofs, some preliminary results are needed. Firstly, we have the following proposition which concerns the multiplicity result for $(\mathcal{SF})$ under the very general assumptions $(f_0)$, $(f_2)$, $(f_3)$ and $(f'_1)$ on $f$, see \cite{Be83-2,Hi10}.

\begin{proposition}
Assume $N\geq2$ and that $f$ satisfies $(f_0)$, $(f_2)$, $(f_3)$ and $(f'_1)$. Then Problem $(\mathcal{SF})$ possesses infinitely many distinct radial solutions $\{v_n\}^{+\infty}_{n=1}$ which satisfy $\|\nabla v_n\|^2_2\to+\infty$ as $n\to+\infty$. Without loss of generality, we may assume that
\begin{equation}\label{equ5.1}
\|\nabla v_n\|^2_2<\|\nabla v_{n+1}\|^2_2~~~~\text{for every}~ n\in\mathbb{N}.
\end{equation}
\end{proposition}

On the other hand, similar to Proposition 2.1 in \cite{Lu15-1}, we have the following proposition.
\begin{proposition}
When $N\geq2$, then $u\in H^1(\mathbb{R}^N)$ is a nontrivial solution to $(\mathcal{KT})$ if and only if there exist a nontrivial solution $v\in H^1(\mathbb{R}^N)$ to $(\mathcal{SF})$ and $t>0$ such that
\begin{equation*}
h(v,t):=M\left(t^{2-N}\|\nabla v\|^2_2\right)t^2=1~~~~\text{and}~~~~u(\cdot)=v(t\cdot).
\end{equation*}
\end{proposition}

Now, under the assumptions of Theorem 1.1, the existence result of infinitely many distinct solutions to $(\mathcal{KT})$ can be proved in a convenient way.

Actually, when $N=2$, let $u_n(\cdot):=v_n(t_n\cdot)$ for every $n\in\mathbb{N}$, where $t_n>0$ is uniquely determined by $h(v_n,t_n)=1$. When $N\geq3$, $(M_1)$ and $(M_3)$ show that, for every $v\not\equiv0$,
\begin{equation*}
h(v,t)\to+\infty~\text{as}~t\to+\infty~~~~\text{and}~~~~h(v,t)\to0^+~\text{as}~t\to0^+.
\end{equation*}
Thus, there exists a positive sequence $\{t_n\}^{+\infty}_{n=1}$ such that $h(v_n,t_n)=1$. In terms of $(M_3)$ and \eqref{equ5.1}, we can also assume that $t^{2-N}_n\|\nabla v_n\|^2_2<t^{2-N}_{n+1}\|\nabla v_{n+1}\|^2_2$ for every $n\in\mathbb{N}$. For such $\{t_n\}^{+\infty}_{n=1}$, set $u_n(\cdot):=v_n(t_n\cdot)$, $n=1,2,\cdots$. From Propositions 5.1 and 5.2 and the fact that $\|\nabla u_n\|^2_2=t^{2-N}_n\|\nabla v_n\|^2_2$, we conclude easily that $\{u_n\}^{+\infty}_{n=1}$ defined as above are the desired solutions for $N\geq2$.

Similarly, under the assumptions of Theorem 1.2, the existence result of finitely many distinct solutions to $(\mathcal{KT})$ can also be proved  from the non-variational point of view. The detailed proof is provided here for reader's convenience.

For every fixed $n\in\mathbb{N}$, let $q\in(0,q_n)$, where
\begin{equation*}
q_n:=\frac{m_0}{1+\underset{1\leq i\leq n}{\max}\lambda\left((2m_0)^{\frac{N-2}{2}}\|\nabla v_i\|^2_2\right)} >0.
\end{equation*}
Obviously, $h\left(v_i,\frac{1}{\sqrt{2m_0}}\right)<1$ for every $i\in\{1,\cdots,n\}$. On the other hand, $(M_1)$ yields that $h(v_i,t)\to+\infty$ as $t\to+\infty$ for every $i\in\{1,\cdots,n\}$. Thus, there exists a positive sequence $\{t_i\}^{n}_{i=1}$ such that $h(v_i,t_i)=1$ for every $i\in\{1,\cdots,n\}$. In terms of \eqref{equ5.1}, we also have that $t^{2-N}_i\|\nabla v_i\|^2_2\neq t^{2-N}_{j}\|\nabla v_{j}\|^2_2$ for every $i,j\in \{1,\cdots,n\}$ and $i\neq j$. For such $\{t_i\}^{n}_{i=1}$, set $u_i(\cdot):=v_i(t_i\cdot)$, $i=1,2,\cdots,n$. Now, it is easy to see  that $\{u_i\}^{n}_{i=1}$ are the desired solutions.

\begin{remark}
In some sense, $J$ can be seen as a suitable perturbation of $I$. Additionally, Proposition 5.2 provides a clear and vital relation between the solutions of $(\mathcal{KT})$ and that of $(\mathcal{SF})$. Thus, in terms of Proposition 5.1, it is natural and well-founded to ask the existence of multiple solutions to $(\mathcal{KT})$.
\end{remark}
\begin{remark}
As we can see in this section, the assumptions on $M$ are mainly used to ensure the existence of $t>0$ such that $h(v,t)=1$. In this procedure, we observe that, when $N\geq3$, the behavior of function $\widehat{M}(t)-\left(1-2/N\right)M(t)t$ at infinity is actually not used, which, in contrast, plays a important role in the variational proofs, see Claim 1 and its proof in Subsection 3.2. This significant difference seems to imply that, in the variational arguments, the boundedness of $\{\|\nabla u_j\|^2_2\}^{+\infty}_{j=1}$ could be established under some weaker assumptions on $M$ or in a more natural way.
\end{remark}

\section*{Acknowledgment}
The author would like to express his sincere gratitude to his advisor Professor Zhi-Qiang Wang for his patient guidance, constant encouragement and timely help. The author also thanks Professor Kazunaga Tanaka for sharing the full text of \cite{Hi10}.


\begin{thebibliography}{100}
\bibitem{Az12} A. Azzollini,
{\it The elliptic Kirchhoff equation in ${\mathbb{R}^N}$ perturbed by a local nonlinearity},
Differ. Integral Equ.  {25} (5-6) (2012), 543--554.
\bibitem{Az15} A. Azzollini,
{\it A note on the elliptic Kirchhoff equation in ${\mathbb{R}^N}$ perturbed by a local nonlinearity},
Commun. Contemp. Math. {17} (2015) 1450039. pages 5.
\bibitem{Az11} A. Azzollini, P. d'Avenia, A. Pomponio,
{\it Multiple critical points for a class of nonlinear functions},
Ann. Mat. Pura. Appl. {190} (2011), 507--523.
\bibitem{Be83-3} H. Berestycki, T. Gallouet, O. Kavian,
{\it Equations de champs scalaires euclidens non lineaires dans le plan},
C. R. Acad. Sci; Paris Ser. I Math. {297} (1983), 307--310.
\bibitem{Be83-1} H. Berestycki, P.L. Lions,
{\it Nonlinear scalar field equations I, Existence of a ground state},
Arch. Rat. Mech. Anal. {82} (1983), 313--346.
\bibitem{Be83-2} H. Berestycki, P.L. Lions,
{\it Nonlinear scalar field equations II, Existence of infinitely many solutions},
Arch. Rat. Mech. Anal. {82} (1983), 347--375.
\bibitem{Fi14} G.M. Figueiredo, N. Ikoma, J.R.S. J\'{u}nior,
{\it Existence and concentration result for the Kirchhoff type equations with general nonlinearities},
Arch. Rat. Mech. Anal. {213} (2014), 931--979.
\bibitem{Hi10} J. Hirata, N. Ikoma, K. Tanaka,
{\it Nonlinear scalar field equations in $\mathbb{R}^N$: mountain pass and symmetric mountain pass approaches},
Top. Meth. Nonlinear Anal. {35} (2) (2010) 253–-276.
\bibitem{Je97} L. Jeanjean,
{\it Existence of solutions with prescribed norm for semilinear elliptic equations},
Nonlinear Anal. {28} (1997), 1633--1659.
\bibitem{Ki83} G. Kirchhoff,
{\it Mechanik},
Teubner, Leipzig, 1883.
\bibitem{Lu15-1} S.-S. Lu,
{\it An autonomous Kirchhoff-type equation with general nonlinearity in $\mathbb{R}^N$}, preprint(2015),
http://arxiv.org/abs/1510.07231.
\end{thebibliography}
\end{document}